\documentclass[preprint]{elsarticle}
\usepackage{lineno,hyperref}
\modulolinenumbers[5]
\usepackage{amssymb}

\makeatletter
\def\ps@pprintTitle{%
  \let\@oddhead\@empty
  \let\@evenhead\@empty
  \def\@oddfoot{\reset@font\hfil\thepage\hfil}
  \let\@evenfoot\@oddfoot
}
\makeatother

\usepackage[ansinew]{inputenc}
\usepackage{amssymb}
\usepackage{amsmath}
\usepackage{amsfonts}
\usepackage{amsthm}

\usepackage{color}

\usepackage{epsfig}

\usepackage{graphics,graphicx}

\usepackage[ruled, lined]{algorithm2e}

\newtheorem{theorem}[equation]{Theorem}

\newtheorem{example}[equation]{Example}

\theoremstyle{definition}
\newtheorem{definition}[equation]{Definition}

\numberwithin{equation}{section}

\begin{document}

\begin{frontmatter}

\title{Derivative of a hypergraph\\ as a tool for linguistic pattern analysis} 

\author[gil,urjc3]{\'Angeles Criado-Alonso}
\ead{angeles.criado@urjc.es}


\author[urjc,dnc,la]{David Aleja}
\ead{david.aleja@urjc.es}

\author[urjc,ccs,dnc,la]{Miguel Romance}
\ead{miguel.romance@urjc.es}

\author[urjc,ccs,dnc,la]{Regino Criado\corref{mycorrespondingauthor}}
\cortext[mycorrespondingauthor]{Corresponding author}
\ead{regino.criado@urjc.es}

\address[gil]{Grupo de Investigaci\'on LIyNMEDIA, Universidad Rey Juan Carlos\\c/Tulip\'an s/n, 28933 M\'ostoles (Madrid)}
\address[urjc3]{Departamento de Filolog\'{\i}a Extranjera, Traducci\'on e Interpretaci\'on, Universidad Rey Juan Carlos, C/Tulip\'an s/n, 28933 M\'ostoles (Madrid), Spain}

\address[urjc]{Departamento de Matem\'atica Aplicada, Ciencia e Ingenier\'{\i}a de los Materiales y Tecnolog\'{\i}a Electr\'onica, Universidad Rey Juan Carlos, C/Tulip\'an s/n, 28933 M\'ostoles (Madrid), Spain}
\address[ccs]{Center for Computational Simulation, Universidad Polit\'ecnica de Madrid, 28223 Pozuelo de Alarc\'on (Madrid), Spain}
\address[dnc]{Data, Complex Networks and Cybersecurity Sciences Technological Institute,  Univ. Rey Juan Carlos, Pza. Manuel Becerra 14, 28028 Madrid, Spain}
\address[la]{Laboratory of Mathematical Computation on Complex Networks and their Applications, Universidad Rey Juan Carlos, Calle Tulip\'an s/n, M\'ostoles, 28933 Madrid, Spain}

\begin{abstract}
The search for linguistic patterns, stylometry and forensic linguistics have in the theory of complex networks, their structures and associated mathematical tools, allies with which to model and analyze texts.  In this paper we present a new model supported by several mathematical structures such as the hypergraphs or the concept of derivative graph to introduce a new methodology able to analyze the mesoscopic relationships between sentences, paragraphs, chapters and texts, focusing not only in a quantitative index but also in a new mathematical structure that will be of singular help to both: detecting the style of an author and determining the language level of a text. In addition, these new mathematical structures may be useful to detect similarity and dissimilarity in texts and, eventually, even plagiarism.
\end{abstract}


\end{frontmatter}

\section{Introduction}\label{sec:intro}


In the last decades the emergence of new structures and models in the field of complex networks and the successive advances in the study and development of their associated tools have made it possible to model the different types of interactions between the diverse parts of a complex system in an efficient and remarkably successful way in practically all areas of knowledge ~\cite{Boca2006,CFGGR,Estrada10,Kivela2013,Newman2010,Wasserman94}. Complex networks have become an essential and indispensable element in the representation of systems for simulating the interactions and relationships between the components of a complex system in domains as diverse as biology, technology, and human social organization~\cite{,Boca2006,Boca2014,CCMR15,Fontoura11,DeA18,DogoMe01,Fe01,Laniru,MR06,Newman2010}.  

It can be said that Network Science can be traced back to the analysis of heterogeneity in real-world complex systems in both nature and function. Thus, the role played by some nodes in these systems is very different from that obtained by the classical Erd\H{o}s-R\'enyi model of random networks, which was a first fundamental milestone in the modeling of real-world complex systems and in the assumption in these models of a first level of heterogeneity \cite{AB2002}. The famous scale-free model made it possible to successfully model real-world complex systems by highlighting the relevant role of nodes with heterogeneous connectivity \cite{AB2002}. A second milestone consisted in the emergence of multilayer network models by taking into account that links could also be heterogeneous in nature \cite{Boca2014}. The third milestone is currently being developed by considering that the heterogeneity of complex systems may affect not only the function of links, but also their nature, since links may be formed by subsets of nodes of different cardinality \cite{Battisetal20}. From collaborative networks to collective social interaction, from trophic networks to biochemical regulatory networks and, in our case, to linguistic networks, many complex systems are produced by considering interactions between more than two nodes simultaneously, making classical network models insufficient. Therefore, the new challenge for the network community is to find new mathematical models that fit multiparty interactions in order to model complex systems with relationships of heterogeneous nature. 

The emergence of new tools that allow large datasets to be handled and analyzed automatically has led to the development of new approaches in many areas of knowledge including text analysis \cite{DeA16,DeA18,Ka15}.

Classical approaches for linguistic analysis of texts were based on simple statistical studies that relied on word frequency \cite{Alt09,DeA18}. 
However, it is noteworthy that in the last decades modern linguistics has received a great advance stemming from the treatment of a language as a system or complex network, having at its disposal in this representation all the tools, measures and procedures to obtain a new, efficient and effective approach to the study of language through complex network that includes qualitative and quantitative aspects~\cite{BorAre10,CL1,DogoMe01,Fe01,Fe05,Liu1,Liu2,MMM16,MLBB16,Soleetal10}.

Therefore, the analysis of linguistic theories supported by the study of specialized corpora and the new approach provided by complex networks makes it possible to obtain certain stylistic and typological characteristics together with some intrinsic properties of languages. The perspective provided by complex networks must go beyond the use of word adjacency or co-occurrence methods which, although they successfully capture the syntactic elements of the texts  \cite{Fe04}, do not have the capacity to represent certain characteristics that develop at the mesoscopic level throughout the text and that have to do with the semantic relationships between the different sentences and paragraphs that compose it.

The linguistic network model we are working with in this manuscript  emerges from the need to work with sentences or paragraphs as a group or collection of certain words in contrast to the type of links considered in previous works where directed and weighted links are used to represent the relationships between linguistic units as in \cite{CL1,MMM16,MR06}. In this work, as  in \cite{CABARC20,CABARC21}, instead of considering the co-occurrence relationship between two adjacent words or linguistic units within a sentence, we will study not only the relationship between sentences (those that share lexical words) but also the relationship between paragraphs or even articles, seeking to characterize, by using network theory parameters, the style of an author or a text as well as the level of language and/or specialization used in the text. This approach leads us to a completely different perspective from the one used, for example, in \cite{DeA18}, where, among many other differences, words are transformed and reduced to their canonical forms and the text is organized in consecutive sets of paragraphs.

In order to apply the tools described in this work, and perform a computer processing on a linguistic corpus understood as a collection of texts collected electronically as a representative sample of texts selected according to certain linguistic criteria \cite{McE}, a corpus of texts composed of 86 extended abstracts (volumes 1-6 of the International Journal of Complex Systems in Science (IJCSS), published between April 2011 and November 2016 (http://www.ij-css.org)) has been considered. This corpus provides us with a total amount of 147637 words as well as 25210 sentences, considered in this study. It should be noted that the unit of analysis from which we start in this work is the sentence, i.e., the words enclosed between two periods \cite{Halli}. In addition, it is important to note that commas and other punctuation marks within the sentence have not been considered for this analysis.

The questions we addressed when we started to write this paper were ``How to characterize the competence level of language used in a text?'' or ``Can the style of an author be determined using specific parameters in the linguistic network under consideration?'' and also, ``What is the combination of words most frequently used in a corpus beyond locating the most relevant individual lexical words?'' or even ``How to determine the most representative words of a text (not necessarily the most frequent)?''

Taking into account that the English language has four main classes of words: nouns, adjectives, verbs, and adverbs, and that other classes of words are prepositions, conjunctions, determiners, interjections, or pronouns, we established in \cite{CABARC20} a four-layer network in order to study a specialty language. In this paper, we will focus on words belonging to the lexical layer, i.e., those significant words (mainly nouns and adjectives) with a specific meaning relevant to the specialty language under study  \cite{CABARC20, CABARC21}.

Therefore in this paper we use the tools and methodology derived from some complex network structures to describe interactions between groups of words, each of these groups being formed by the lexical words belonging to a specific sentence in the analyzed corpus (syntagmatic approach, from the Greek ``$\sigma\iota \upsilon\tau\alpha\gamma\mu\alpha$", syntagma: ``assembled group"). It is therefore important to note that the syntagmatic approach, which corresponds to the analysis presented in this paper, is different from the paradigmatic approach used in other works of computational linguistics \cite{CL1}. 

Since the syntagmatic relationship is based on the interrelationships of words in a linguistic structure \cite{CABARC20,CABARC21,Sinclair91}, it makes sense to consider the relationships between sets of two, three or more significant words that appear in the same sentence, paragraph, abstract or article and that in some way characterize a text as belonging to an author, or discriminate the level of language used in it, as well as those other words and relations that allow distinguishing it from texts belonging to other authors or that use a different level of language.  

The methodology presented here makes it possible to determine the level of language used in a text as well as the style of an author and also to analyze and order sentences, abstracts, paragraphs and texts (sets of words) according to their importance, having mind their interrelationships in the context of the multilayer network structure defined in \cite{CABARC20,CABARC21} as well as to extract new features of a text from the relationships between significant sets of words in the text.

High-order networks or hypergraphs are the natural generalization of networks that takes into account the fact that a link can connect more than two nodes. Interest in this type of network is growing due to the inability of classical graphical representations to describe group interactions. Their applicability goes beyond the field of social sciences \cite{Benson19,Lambi19,Torresetal21} and the study of group interactions, public cooperation or opinion formation. In our case, we will consider its applicability in the field of linguistics and specialty languages beyond other approaches based on classical complex networks, multiplex networks or multilayer networks \cite{BorAre10,CL1,Fe01,Fe05,Liu1,Liu14,Liu2,MR06,Sole05,Zipf}.

As it can be easily understood, a property referring to a finite set of objects (in our case, the nodes of a network), is completely characterized by the subset of elements that satisfy it, which in this case will be represented by the hyperedge formed by these elements, so that it will be possible to compare and relate properties of the nodes and the network by studying and analyzing the corresponding hypergraph.

Thus, studying the relationships between the properties of the nodes consists of mathematically analyzing the properties and typical parameters of the associated hypergraph.  Therefore, the applications of this methodology to the field of linguistics range from the characterization of an author's style to the detection of plagiarism, including the detection and identification of the same concept expressed in a different way. To this end, starting, in the first instance, from the identification of a sentence of our corpus with the hyperedge formed by the set of lexical words of that sentence, the hypergraph will be constructed in which the nodes will be all the lexical words of the corpus and the hyperedges all the sentences of the corpus, defining the concept of derivative of two words with respect to a set of hyperedges and the degree of independence of two words of a text with respect to that set of hyperedges.
This study can be extended in more depth by considering as hyperedges, successively, the sets of nodes formed by the lexical words of a paragraph, an abstract or even a chapter, taking the corresponding sequence of parameters as a feature of the text and pointing to new applications of this structure.

The structure of the paper is as follows.  After this introduction, in Section~\ref{sec:basic} some basic concepts and a summary of the most important relationships between the line graph the dual hypergraph, the bipartite graph associated to a certain hypergraph and its corresponding matrices are introduced. Section~\ref{sec:deri}  is devoted to introduce the concept of derivative of a hypergraph with respect to a set of nodes and to establish the definition of the homogeneity graph of a hypergraph obtaining some remarkable results related to this new structure.  In Section~\ref {sec:multi} we apply the mathematical concepts and structures defined in the previous sections to obtain tools to characterize the style and level of a text belonging to the linguistic hypergraph considered. In Section~\ref{sec:experiments} the lexical density of the set of texts that make up the analyzed corpus is studied, and some numerical experiments and computational results are presented by using three different algorithms to illustrate the diverse types of relationships that can be established between sentences within a text and their relative importance. Section~\ref{sec:seeking} is devoted to apply the instruments and tools developed in order to obtain distinctive characteristics that allow us to distinguish the styles of the different authors and linguistic competence levels of the written texts included in the corpus considered. Finally in Section~\ref{sec:conclu} we present some conclusions of this work. 

\section{Basic concepts and some preliminary results}\label{sec:basic}
A network (or graph) $G=(X,E)$ is just a finite set of vertices (or nodes) $X=\{1,...,N\}$ connected by a set of edges (or links between certain pairs of nodes) $E=\{e_1,\cdots,e_m\}$. If the edges have a direction, we will say that $G$ is a directed network (or digraph). In the sequel, we will denote by $e_{ij} \in E$  the link between the nodes $i$ and $j$, although sometimes we will also denote the edge $e_{ij}$  by $\{i,j\}$ or, if $G$ is a directed network, by $i \rightarrow j$. Finally, a weighted network is a graph in which each edge $e_{ij}$  has an associated numerical value $w(e_{ij} )=w_{ij}$ called its weight. In the same way, following  \cite{Berge1}, a hypergraph $\mathcal{H}=(X,\varepsilon)$ is a finite set of vertices (or nodes) $X=\{1,...,N\}$ and a collection $\varepsilon =\{h_1,h_2,\ldots,h_n\}$ of subsets of $X$ such that $h_i\neq\emptyset \quad (i=1,2,\ldots,n)$ and $X=\bigcup_{i=1}^{n} h_i.$ Each of these subsets is called a hyperedge. In this way,  hypergraphs appeared as the natural extensions of graphs to describe group interactions.  In the following sections, the study is developed with undirected graphs and hypergraphs, though some of the definitions can be easily extended to the directed case.

In order to carry out our study it is necessary to introduce the concepts of linegraph and dual hypergraph of a hypergraph. In this regard it should be noted that the concept of linegraph $L(G)$  associated to a graph $G=(X,E)$  was introduced by H. Whitney in 1932 \cite{Whitney32} and extended for higher order networks
by J.C. Bermond et al. in 1977  \cite{Bermondetal77,Tyetal98}. It is important to point out that the study of these structures, as well as the relationships between them and their applications, has been increasing steadily in recent years (see, for example, \cite{bagga04,Benson19,CFGR14,CFGRBM16,EL1,EL2,ran19}).

So, if $\mathcal{H}=(X,\varepsilon)$ is a hypergraph, the linegraph associated to $\mathcal{H}$ is the graph $L(\mathcal{H})=(\varepsilon, E')$, where if $h_i,h_j \in \varepsilon$, then 
\[
\{h_i,h_j\} \in E' \Leftrightarrow h_i \cap h_j \neq \emptyset.
\]
It is also notorious that the linegraph $L(\mathcal{H})$ of a hypergraph $\mathcal{H}$  is a graph even though $\mathcal{H}$  is a hypergraph. Note that this concept is a particular case of the concept of intersection graph \cite{ran19}. On the other hand, it is also possible to consider the dual hypergraph of a hypergraph: 
if  $\mathcal{H}=(X,\varepsilon)$ is a hypergraph, the dual hypergraph associated with $\mathcal{H}$ is the hypergraph $\mathcal{H} ^*=(\varepsilon,X')$ in such a way that if  $X=\{1,...,N\}$, then  $X'=\{v_1,...,v_N\}$ where $v_i=\{h_j| i\in h_j\}$, $i=1,...,N$. It is not difficult to verify that  $(\mathcal{H}^{*})^*= \mathcal{H}$. Moreover, if $I$ is the incidence matrix of $\mathcal{H}$, then its transpose matrix $I^t$ is the incidence matrix of  $\mathcal{H}^*$. In this context, to concretize the relationship between $L(\mathcal{H})$ and $\mathcal{H}^*$, we consider the function $\Pi_2$  that turns a hypergraph  $\mathcal{H}=(X,\varepsilon)$ into a graph  $\Pi_2(\mathcal{H})=(X,E')$ as follows:  
\[
\{i,j\} \in E' \Leftrightarrow \exists h \in \varepsilon \mid ~ i,j \in h.
\]

So, for any hypergraph $\mathcal{H}$ we have that  $\Pi_2(\mathcal{H}^*)=L(\mathcal{H})$. Furthermore, if  $G=(X,E)$ is a graph, with $X=\{1,...,N\}$,  we can also consider the dual hypergraph $G^*=(E,\varepsilon)$ of $G$ where $\varepsilon=\{h_1,.... h_n\}$ and $\forall i \in \{1,...n\}$ we consider the corresponding hyperedge $h_i=\{e_j \in E |~ i \in e_j\}$, and also  $\Pi_2(G^*)=L(G)$.

Now, if we denote by  $I(\mathcal{H})$ the incidence matrix of $\mathcal{H}$, then it is not difficult to verify that
\[
I(\mathcal{H})^t\cdot I(\mathcal{H})=\widetilde{A(\mathcal{H})}=(\widetilde{a_{ij}})\in \mathbb{R}^{|\varepsilon| \times |\varepsilon |}
\]
and
\[
I(\mathcal{H})\cdot I(\mathcal{H})^t=A(\mathcal{H})=(a_{ij})\in \mathbb{R}^{N \times N},
\]
where
\[\widetilde{a_{ij}}=\left\{
\begin{array}{cl}
|h_i| & {\rm{if}}\quad i=j, \\
|h_i \cap h_j| & {\rm{if}}\quad i \ne j,
\end{array}
\right.
\]
and
\begin{equation}\label{eq:FreqMatrix}
a_{ij}=\left\{
\begin{array}{cl}
|\{h \in \varepsilon  \mid ~i \in h\}| & {\rm{if}}\quad i=j, \\
|\{h \in \varepsilon  \mid ~ i, j \in h\}| & {\rm{if}}\quad i \ne j.
\end{array}
\right.
\end{equation}

In fact, if  we consider in addition the bipartite network $B(\mathcal{H})$ associated to the  hypergraph $\mathcal{H}=(X,\varepsilon)$ defined by $B(\mathcal{H})=(X\cup \varepsilon, E(\mathcal{H}))$ then its adjacency matrix is given by
\[
A_{B(\mathcal{H})}= \left(
\begin{array}{c|c}
0 & I(\mathcal{H}) \\
\hline I(\mathcal{H})^t& 0
\end{array}
\right)
\]
and 
\[
(A_{B(\mathcal{H})})^2= \left(
\begin{array}{c|c}
A(\mathcal{H}) & 0 \\
\hline 0& \widetilde{A(\mathcal{H})}
\end{array}
\right).
\]

The matrix $A(\mathcal{H})=(a_{ij})$  is called the frequency matrix of relations between the elements (nodes) of the hypergraph $\mathcal{H}$  (see \cite{Gorba88}).


\section{Hypergraphs and Derivative graph}\label{sec:deri}
Quantifying the similarity between two models or structures is one of the most important aspects that has contributed to the development of theories and models in science and technology.
There are multiple works whose objective is to model generic data sets in the field of complex networks in order to, by using the constructed model, study the level of similarity or coincidence of such data \cite{Brus,LFC22JC,Vij}. 
 Thus, since the introduction of Jaccard's index in 1901 \cite{Jacc}, through different adaptations and generalizations of this concept  \cite{LFC21,Vij}, several types of indexes  and generalizations have been established with the aim of quantifying the similarity between two sets or mathematical structures \cite{Brus,Ham,Vij,LFC21,LFC22PA,LFC22JC}.

The basic Jaccard index to compare the degree of coincidence or similarity between two sets $A$ and $B$ can be obtained from the formula 
\[
\mathcal{J}=\frac{|A\cap B|}{|A\cup B|}.
\]  
The different applications of the Jaccard index along time made possible the development of new indexes, improving the accuracy of the original results. So, the overlap index and  the coincidence similarity \cite{LFC22, LFC22PA,LFC22JC,Vij} are examples of additional indexes that allow to establish similarity between certain types of models and structures, including approaches aimed at quantifying similarity between paragraph contents using the concept of multisets \cite{LFC22}.
 
In our case, we are going to introduce a methodology to analyze and quantify the similarity between two nodes $i,j$ of a hypergraph, applying it to the study of the linguistic network built through the corpus under study.  

In this section we are going to introduce the concept of derivative  graph of a hypergraph with the idea of associating not only a numerical index that allows us to quantify the heterogeneity and absence of similarity between the corresponding hyperedges, but also to associate a structure (in this case a graph) to characterize the heterogeneity and dissimilarity of the elements of the hypergraph under consideration. Now, we are in a good position to establish the concept of derivative graph of a hypergraph over a pair of nodes:

\begin{definition}\label{deri}
Given  a hypergraph $\mathcal{H}=(X,\varepsilon)$, with $A(\mathcal{H})=(a_{ij})\in \mathbb{R}^{N \times N}$, we will call the derivative hypergraph of  $\mathcal{H}$ with respect to the pair of nodes $i,j \in X$ as the numerical value  $\frac{\partial \mathcal{H}}{\partial \{i,j\}}$ obtained by applying the following formula
\begin{equation}\label{eq:deri}
\frac{\partial \mathcal{H}}{\partial \{i,j\}}= \frac{a_i-a_{ij}+a_j-a_{ij}}{a_{ij}}= \frac{a_i-2a_{ij}+a_j}{a_{ij}}.
\end{equation}
\end{definition}

Obviously, if there is not a hyperedge $h \in \varepsilon $ such that $i,j \in h$, we will have 
$\frac{\partial \mathcal{H}}{\partial \{i,j\}}= \infty$, and if $\forall h \in  \varepsilon ~~(i \in h \Leftrightarrow j \in h)$ then we will have $\frac{\partial \mathcal{H}}{\partial \{i,j\}}=0$.
Note that $\forall i,j \in X$ we have that  $\frac{\partial \mathcal{H}}{\partial \{i,j\}} \ge 0$.

It is important to point out that the above definitions can be extended without difficulty to the context of a collection of sets (which would play the role of the hyperedges) and of the elements (respectively the nodes) of the sets of that collection.

If we now consider each hyperedge $h \in \varepsilon$  as a property or a feature that a node may or may not have, or even as an event or affair in which a particular node may or may not participate, so that the entire hypergraph is a set of features or events, the value of  $\frac{\partial \mathcal{H}}{\partial \{i,j\}}$ characterizes the (relative) heterogeneity of the properties $\varepsilon$ satisfied simultaneously by nodes $i$ and $j$, or the intensity of participation  of the nodes $i$ and $j$ in the set of events $\varepsilon$. Moreover, the smaller the value of the derivative of the network with respect to the set of events over the pair of nodes $i,j$, the greater  identification and similarity between the corresponding nodes $i,j$ with respect to the considered set of events (in fact, if $\frac{\partial \mathcal{H}}{\partial \{i,j\}}=0$,  these nodes, which participate in exactly the same hyperedges, are so similar that they are, from the point of view of  $\mathcal{H}$ indistinguishable). In other words, the higher the value of the derivative is, the greater the degree of unequal participation of the nodes in the hyperedges. Thus, it makes sense to give the following definition:

\begin{definition}
Given  a hypergraph $\mathcal{H}=(X,\varepsilon)$ and $i,j \in X$, we will call  degree of independence of $i$ and $j$ with respect to $\mathcal{H}$  the numerical value of  $\frac{\partial \mathcal{H}}{\partial \{i,j\}}$.
\end{definition}

\begin{definition}
Given  a hypergraph $\mathcal{H}=(X,\varepsilon)$, the derivative graph $\partial \mathcal{H}$  of  $\mathcal{H}$ is the weighted graph obtained by considering the derivative of   $\mathcal{H}$ with respect all the pairs of nodes $i,j \in X$, and by setting $\forall i,j \in X$ the corresponding numerical value of $\frac{\partial \mathcal{H}}{\partial \{i,j\}}$ on the edge $\{i,j\}$,  in such a way that if $\frac{\partial \mathcal{H}}{\partial \{i,j\}}= 0$, then the nodes $i$ and $j$ collapse into a single node $(ij)$, and having in mind that if  $\frac{\partial \mathcal{H}}{\partial \{i,j\}}= \infty$, then the edge $\{i,j\}$ does not exist in the derivative graph.
\end{definition}

Globally, it can be said that the derivative graph  $\partial \mathcal{H}$ gives us a representation of  the degree of heterogeneity of participation of nodes on the different hyperedges of $\mathcal{H}$. 

Assuming that if $k$ is any positive number then $\frac{k}{0} =+ \infty$ and  $\frac{k}{\infty} =0$, for continuity and consistency sake of the established concepts, we are interested in defining the homogeneity matrix and homogeneity graph of a hypergraph:

\begin{definition}\label{hom}
Given  a hypergraph $\mathcal{H}=(X,\varepsilon)$, we will call  homogeneity matrix of  $\mathcal{H}$,  to the matrix  $H(\mathcal{H})=(h_{ij})\in \mathbb{R}^{N \times N}$ defined by 
\[
h_{ij}=\left\{
\begin{array}{cl}
0 & {\rm{if}}\quad i=j, \\
\displaystyle\frac{1}{\frac{\partial \mathcal{H}}{\partial \{i,j\}}} & {\rm{if}}\quad i \ne j.
\end{array}
\right.
\]
\end{definition}

\begin{definition}
Given  a hypergraph $\mathcal{H}=(X,\varepsilon)$, the homogeneity graph $HG( \mathcal{H})$  of  $\mathcal{H}$ is the weighted graph with the same nodes and edges as  $\partial \mathcal{H}$ , but considering as the weight of each edge the inverse value of the weight corresponding to the derived graph  $\partial \mathcal{H}$.
\end{definition}

At this point it is remarkable that the application of the PageRank algorithm on the homogeneity graph $HG( \mathcal{H})$ will allow us to extract the most representative nodes of the hypergraph, in the sense that the nodes located in the first places of the ranking obtained will be the ``most similar'' (in the sense that underlies the definition of homogeneity graph) to each other and to the rest of the nodes of the hypergraph as it will be shown in Section~\ref{sec:experiments}.

\begin{figure}[ht]
  \begin{center}
    \includegraphics[width=1\textwidth]{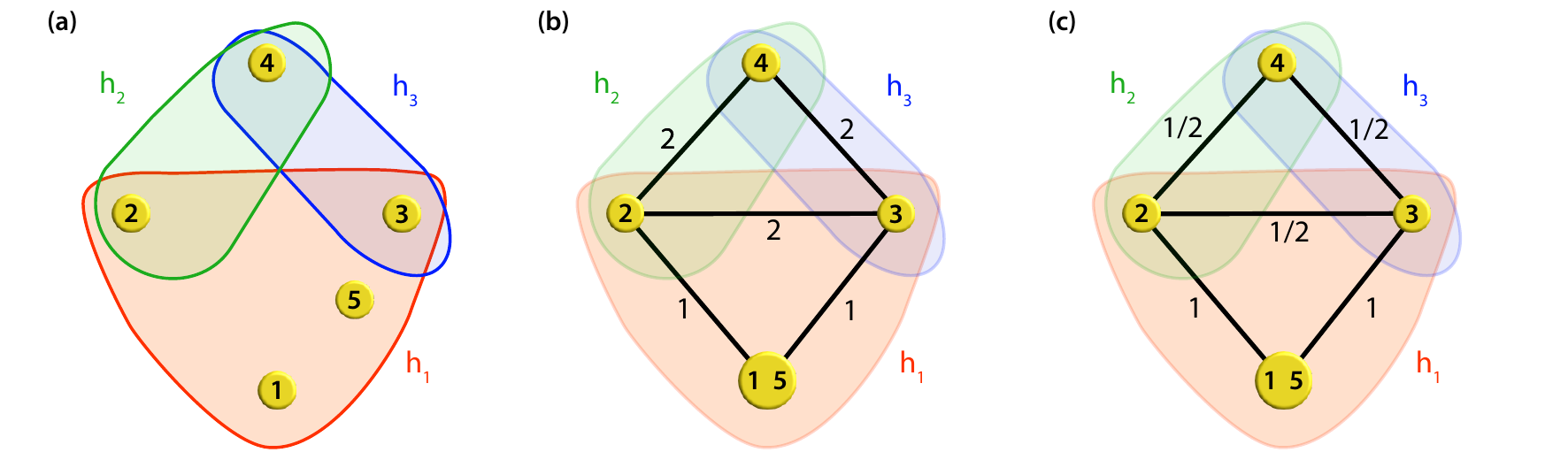}
  \end{center}
 \caption{Hypergraph $\mathcal{H}$ (panel (a)), its derivative graph $\partial \mathcal{H}$ (panel (b)) and  its homogeneity graph $HG( \mathcal{H})$ (panel (c)).}
\label{Fig5}
 \end{figure}

To clarify the concepts and ideas introduced, let's examine the following example:

\begin{example}
Consider the hypergraph  $\mathcal{H}=(X,\varepsilon)$, where $X=\{1,2,3,4,5\}$, $\varepsilon=\{h_1, h_2, h_3\}$, and  $h_1=\{1,2,3,5\},\,h_2=\{2,4\},\,h_3=\{3,4\}$, represented in panel (a) of  Figure~\ref{Fig5}. We have that 

\[
 I(\mathcal{H})^t=\left(\begin{array}{ccccc}
1&1&1&0&1\\
0&1&0&1&0\\
0&0&1&1&0
\end{array}\right),
\]

\[
I(\mathcal{H})\cdot I(\mathcal{H})^t=A(\mathcal{H})=\left(\begin{array}{ccccc}
1&1&1&0&1\\
1&2&1&1&1\\
1&1&2&1&1\\
0&1&1&2&0\\
1&1&1&0&1
\end{array}\right).
\]

 The values of the derivatives of  $\mathcal{H}$ with respect  to all the pair of nodes of $G$ are, respectively:

\[\
\begin{split}
\frac{\partial \mathcal{H}}{\partial \{1,5\}}&=0,\, 
\frac{\partial \mathcal{H}}{\partial \{1,2\}}=1,\,
\frac{\partial \mathcal{H}}{\partial \{1,3\}}=1,\,
\frac{\partial \mathcal{H}}{\partial \{1,4\}}=+\infty,\,
\frac{\partial \mathcal{H}}{\partial \{2,3\}}=2,\\
\frac{\partial \mathcal{H}}{\partial \{2,4\}}&=2,\,
\frac{\partial \mathcal{H}}{\partial \{2,5\}}=1,\,
\frac{\partial \mathcal{H}}{\partial \{3,4\}}=2,\,
\frac{\partial \mathcal{H}}{\partial \{3,5\}}=1,\,
\frac{\partial \mathcal{H}}{\partial \{4,5\}}=+\infty.
\end{split}
\]
so that the derivative graph $\partial \mathcal{H}$ is the one represented in part (b) of  Figure~\ref{Fig5} and the  homogeneity matrix of  $\mathcal{H}$  is:

\[
H(\mathcal{H})=\left(\begin{array}{ccccc}
0&1&1&0&\infty\\
1&0&1/2&1/2&1\\
1&1/2&0&1/2&1\\
0&1/2&1/2&0&0\\
\infty&1&1&0&0
\end{array}\right).
\]

Note that the edge $\{1,4\}\in E$  has been removed in the derivative network $\partial \mathcal{H}$ and that nodes $1$ and $5$ have collapsed into a single node in the obtained network. So, the adjacency matrix of the homogeneity graph $HG(\mathcal{H})$  is:

\[
\left(\begin{array}{ccccc}
0&1&1&0\\
1&0&1/2&1/2\\
1&1/2&0&1/2\\
0&1/2&1/2&0\\
\end{array}\right),
\]
where the set of nodes of $HG(\mathcal{H})$  is $(\{(1,5),2,3,4\}$ ordered as they appear (panel (c) of Figure~\ref{Fig5}).
\end{example}

Thus, in panel (a) of  the Figure~\ref{Fig5} it can be observed the original hypergraph $\mathcal{H}$, in part (b) its  derivative graph  $\partial \mathcal{H}$ and in panel (c)  its corresponding homogeneity graph $HG( \mathcal{H})$.

It is worth noting that, in a similar way as it has been done in Definition~\ref{deri},  it is possible to establish the derivative of a hypergraph with respect to a set of three or more nodes as follows:
\[
\frac{\partial \mathcal{H}}{\partial \{i,j,k\} } = \frac {1}{a_{ijk}} \cdot \left( \sum_{\begin{subarray}{l} r \in \{i,j,k\} \\\end{subarray}} \!\!\!\!\!a_r     \enspace-\enspace2\!\!\!\!\!  \sum_{\begin{subarray}{l} r,s \in \{i,j,k\} \\ r \ne s \end{subarray}} \!\!\!\!\!a_{rs}  + 3 a_{ijk} \right),
\]
where $a_{ijk}=|\{h \in \varepsilon|\, i, j,k \in h\}|$, and the same type of formula can be obtained for sets of nodes of higher cardinality.

Note that the same idea can be extended to the definition of degree of independence of several nodes as follows:  Given a hypergraph $\mathcal{H}=(X,\varepsilon)$, and  $i_1,...i_n \in X$,  the  degree of independence of $i_1,...,i_n$ in $\mathcal{H}$ is the numerical value
$\frac{\partial \mathcal{H}}{\partial  \{i_1,...i_n \}}$.

Finally, it is remarkable that the use of the PageRank algorithm on the homogeneity graph will allow us to extract a ranking of the most representative individuals (or nodes) of either the hypergraph or the network under consideration.

To conclude this section, it must be noted that when both graphs and hypergraphs are used simultaneously to model certain complex systems, it is sometimes very useful to analyze how these structures interact and overlap using the tools introduced in this section.
 In this regard, it should be noted that the tools introduced in this section can be used to capture intrinsic and mesoscopic characteristics of a graph and to define new invariants of graphs and isomorphic networks. For example, given a graph $G=(X,E)$, we can consider the hypergraph $\mathcal{H}=(X,\varepsilon)$ such that each of its hyperedges is formed by all the nodes that are part of a cycle, or by all the nodes that are part of a spanning tree of $G$.   
The most accurate framework to work with the overlapping of these structures is the use of hyperstructures.

In  \cite{CRV10}  we can find a first definition of the concept of hyperstructure  as follows:

\begin{definition}[\cite{CRV10}] 
Given a graph $G=(X,E)$ with $N$ vertices and $m$ edges and a hypergraph $\mathcal{H}=(X,\varepsilon)$,  a \emph{hyperstructure} $S=(X,E,\mathcal{H})$ is a triple formed
by the vertex set $X$, the edge set $E$ and the hyperedge set $\mathcal{H}$. The hyperstructure $S$ is said to be compatible if  for every edge $e=\{v,w\} \in E$ there exists a hyperedge $h\in \varepsilon$ such that $v,w \in h$.
\end{definition}

It is not difficult to prove the following result:

\begin{theorem}\label{sec:theorem1}
Let $S=(X,E,\mathcal{H})$ be a hyperstructure, $L(G)=(E, E')$ the linegraph of $G=(X,E)$ and $\Pi_2(\mathcal{H})=(X,E'')$.  If $S$ is compatible, then $S'=(E,E',\mathcal{H})$ and $S''=(X,E'',\mathcal{H})$  are also hyperstructures. 
\end{theorem}

It is important to highlight that by using the idea of derivative we have introduced in this paper we can examine and determine the uniformity of participation of two, three or more nodes in the considered structure or hyperstructure, or even the binary relationships (edges) between participants of a certain event by simply considering a suitable hyperstructure in which the nodes be the edges of the original graph under consideration.

Now, we can define the derivative graph of a weighted hyperstructure:
\begin{definition}
Given a hyperstructure $S=(X,E,\mathcal{H})$, where $G=(X,E,W)$  is a weighted graph and $\mathcal{H}=(X,\varepsilon)$, if $w_{ij}$ denotes the weight of the edge $e=\{i,j\} \in E$, then we will call the derivative of $e$ with respect to the  hyperstructure $S$ the numerical value obtained by applying the following formula
\[
\frac{\partial e}{\partial S} = w_{ij} \cdot \left(\frac{a_i-2a_{ij}+a_j}{a_{ij}}\right).
\]
\end{definition}

Obviously, if there is not a hyperedge $h \in \mathcal{H} $ such that $e=\{i,j \} \in h$, we will have $\frac{\partial e}{\partial S}= \infty$.  On the other hand, it is evident that if a hyperstructure is compatible, the derivative of any edge with respect to $S$ cannot be equal to $+\infty$.

\begin{definition} \label{hip}
Given  a hyperstructure $S=(X,E,\mathcal{H})$, where $G=(X,E,W)$  is a weighted graph and $\mathcal{H}=(X,\varepsilon)$, if $w_{ij}$ denotes the weight of the edge $e=\{i,j\}$, then the derivative graph of $G$ with respect to  $S$ is the weighted graph  $\frac{\partial G}{\partial S}$ obtained by setting $\forall e \in E$ the corresponding numerical value of  $\frac{\partial e}{\partial S}$ on the edge $e=\{i,j\}$, in such a way that if $\frac{\partial e}{\partial S}= 0$, then the nodes $i$ and $j$ collapse into a single node $(ij)$.
\end{definition}

\begin{figure}[ht]
  \begin{center}
    \includegraphics[width=0.7\textwidth]{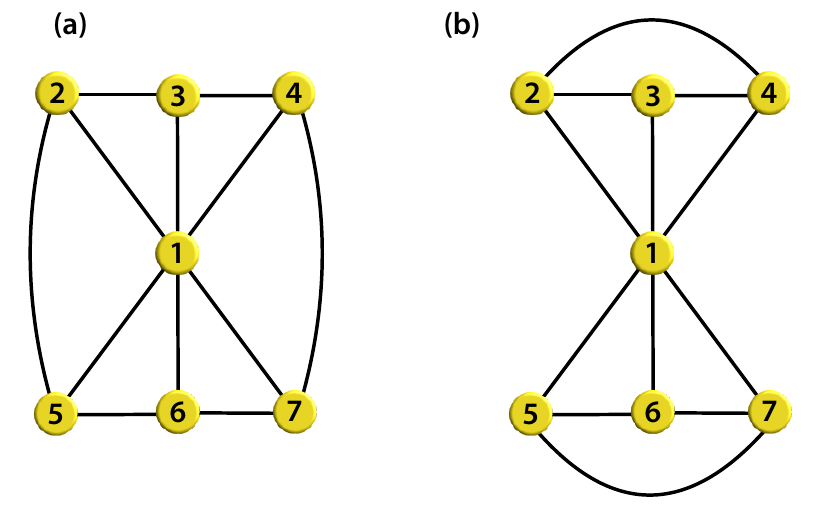}
  \end{center}
 \caption{Two graphs with the same degree distribution: Graph $G$ (panel (a)) and graph $G'$ (panel (b))}
\label{Fig7}
 \end{figure}

As a direct application of the definition, note that if we consider the graphs $G=(X,E)$  (panel (a) of Figure~\ref{Fig7}) and $G'=(X,E')$ (panel (b) of Figure~\ref{Fig7}) and  the hyperstructures  $S=(X,E,\mathcal{H})$ and $S=(X,E',\mathcal{H'})$ such that each of their hyperedges  is composed by all the nodes belonging to a cycle formed by three or more nodes of $G$ and $G'$ respectively, then the derived graphs  $\frac{\partial G}{\partial S}$ and  $\frac{\partial G'}{\partial S'}$ are completely different since, for example, 
\[
\frac{\partial \{1,2\}}{\partial S'}=\frac{20}{8}=\frac{5}{2}.
\]
On the other hand, as can be seen, 
\[
\frac{\partial \{1,2\}}{\partial S}=\frac{62}{20}=\frac{31}{10},
\]
and, obviously,
\[
\frac{5}{2} \ne \frac{31}{10}.
\]

Note that Definition~\ref{hip} allows us to iterate the derivatives with respect to a hyperstructure, because if the graph derived from the hyperstructure is $\frac{\partial G}{\partial S}=(X',E',W')$ y $S'=(X',E',\mathcal{K})$, then we can consider  the mixed derivatives of a graph $G$ with respect to two different hyperstructures (which may eventually be the same) $S$ and $S'$ (in this order) as 
\[ 
\frac{\partial^2 G}{\partial S' \partial S} = \frac{\partial}{\partial S'} \left( \frac{\partial G}{\partial S}\right).
\] 
It is obvious that the successive derivative graphs obtained by deriving respect a suitable chain of two or more  hyperstructures allow to obtain characteristics and properties of the system or model under study related to the absence of similarity between the nodes.







\section{A linguistic hyperstructure based on the  lexical layer within a multilayer linguistic network model} \label{sec:multi}
We are now ready to show the potential applications of  the defined mathematical structures and tools to the linguistic analysis of texts, looking for the identification of signs and specific features of a style or competence level of language considering the most significant words and their relationships. 
It can be said that the English language has four major word grammar categories: nouns, adjectives, verbs, and adverbs. Other word classes are prepositions, conjunctions, determiners, interjections or pronouns \cite{Huddleston}. On this basis described in \cite{CABARC20,CABARC21} we have built a methodology close to supervised machine learning consisting of dividing the words of the corpus under study into a multilayer network  \cite{Boca2014} composed by four layers:  lexical layer, verb layer, linking layer and remaining words layer.

In order to discriminate between the terms (words) and to assign them to one or another layer, a completely lexical linguistic decision was made according to the criteria of several experts. Thus, the terms (words) of the corpus have been distributed in the different layers according to their morphological and  lexical properties. Some other linguistic aspects, such as the specific terminology of a specialty language and the different combinations of words that give rise to new meanings (called ``linguistic collocations") have also been successfully studied and  modeled in \cite{CABARC20,CABARC21}. 

In the model established in \cite{CABARC20} interlayer relations are the basic grammatical relations in a sentence, for example, the interaction between layers that facilitates the formation and description of specialty verbs (e.g. ``cluster together''). On the other hand, throughout the present work, we will consider the sentences as the unit under study, identifying each sentence in the corpus (set of words located between two periods) with the subset of lexical words appearing in that sentence. 

For this reason, throughout this work we are going to focus on the words (nodes) located in the lexical layer. At this point, it is remarkable that in the lexical layer many words can act as verbs when we analyze texts written by authors with higher language skills. For example, within the sentence ``model a network'', the word ``model'' is a verb, but in the expression ``network model'' the term ``model'' is a noun. 

In order to set our approach,  the model  of the corpus analyzed is considered as a set of texts formed by sentences (set of lexical words between two periods). In fact, from a practical and computational point of view, each sentence is identified with the set of lexical words that compose it.  This way, let us consider the hyperstructure in which the nodes are the lexical words, the edges between these nodes are established when these words appear in the same sentence, and the set of hyperedges is the set of sentences that constitute the corpus. 

It is important to point out that the linguistic hyperstructure considered is a compatible hyperstructure, since the edges are established between words that appear in the same sentence. Therefore, from the Theorem~\ref{sec:theorem1} it is possible to study both the hyperstructure in which the nodes are the words and the hyperedges are the sentences and, in a complementary way, the hyperstructure in which the nodes are the edges between words (dual graph of the original graph) and the hyperedges are also the sentences. 

On the other hand, by considering paragraphs as a set of sentences, and the extended abstracts of our corpus as a set  of paragraphs, we can add to this model new linguistic hyperstructures that undoubtedly allow us to characterize a text or set of texts from the derivatives of the corresponding graphs and hypergraphs respectively.

In order to illustrate how useful are the tools presented in the context of the linguistic analysis of texts, let us consider a text in which the same sentence is repeated over and over again. In that case, by deriving the linguistic hypergraph formed by the set of all the repeated sentences with respect to the lexical words of the sentence repeated over and over in all those sentences, the derivative graph will collapse to a single node.

So, by calculating the derivative graph from the linguistic hypergraph composed by all the sentences of a corpus or a text,  we will obtain the degree of similarity between the sentences of that text, and also the greater or lesser degree of difference between all the sentences forming such text (or corpus), with the peculiarity that these quantitative measures are represented in the corresponding derivative graph.
 
Consequently, the derivative graph of a text or a set of texts is  a quantitative and qualitative structure of such text that is a specific feature of that text (or set of texts) for real, which may be considered, in certain cases, like a signature or specific characteristic of the style of an author.

When analyzing the hypergraph $\mathcal{H}$ formed by all the sentences of the corpus under study, we obtained $127$ pairs of words that appear in exactly the same sentences. Thus, for example  
\[
\begin{split}
\frac{\partial \mathcal{H}}{\partial \{monte,carlo\}}&=
\frac{\partial \mathcal{H}}{\partial \{differential,rungekuttta\}}\\
&=\frac{\partial \mathcal{H}}{\partial \{oscillatory, asynchronous\}}=0.
\end{split}
\] 
It is important to note at this point that, if three or more words in the corpus analyzed appear in exactly the same sentences, these words have collapsed into a single node. This has happened in $13$ cases.

Finally, and by way of illustrative example, we will point out that
\[
\begin{split}
\frac{\partial \mathcal{H}}{\partial \{network,system\}}&=16.33,\\ 
\frac{\partial \mathcal{H}}{\partial \{language,formal\}}&=\frac{\partial \mathcal{H}}{\partial \{connectance,asymmetry\}}=2,\\
\frac{\partial \mathcal{H}}{\partial \{process,graph\}}&=28.14,\\
\frac{\partial \mathcal{H}}{\partial \{placing,model\}}&=+\infty.
\end{split}
\]

\begin{figure}[ht]
  \begin{center}
    \includegraphics[width=1\textwidth]{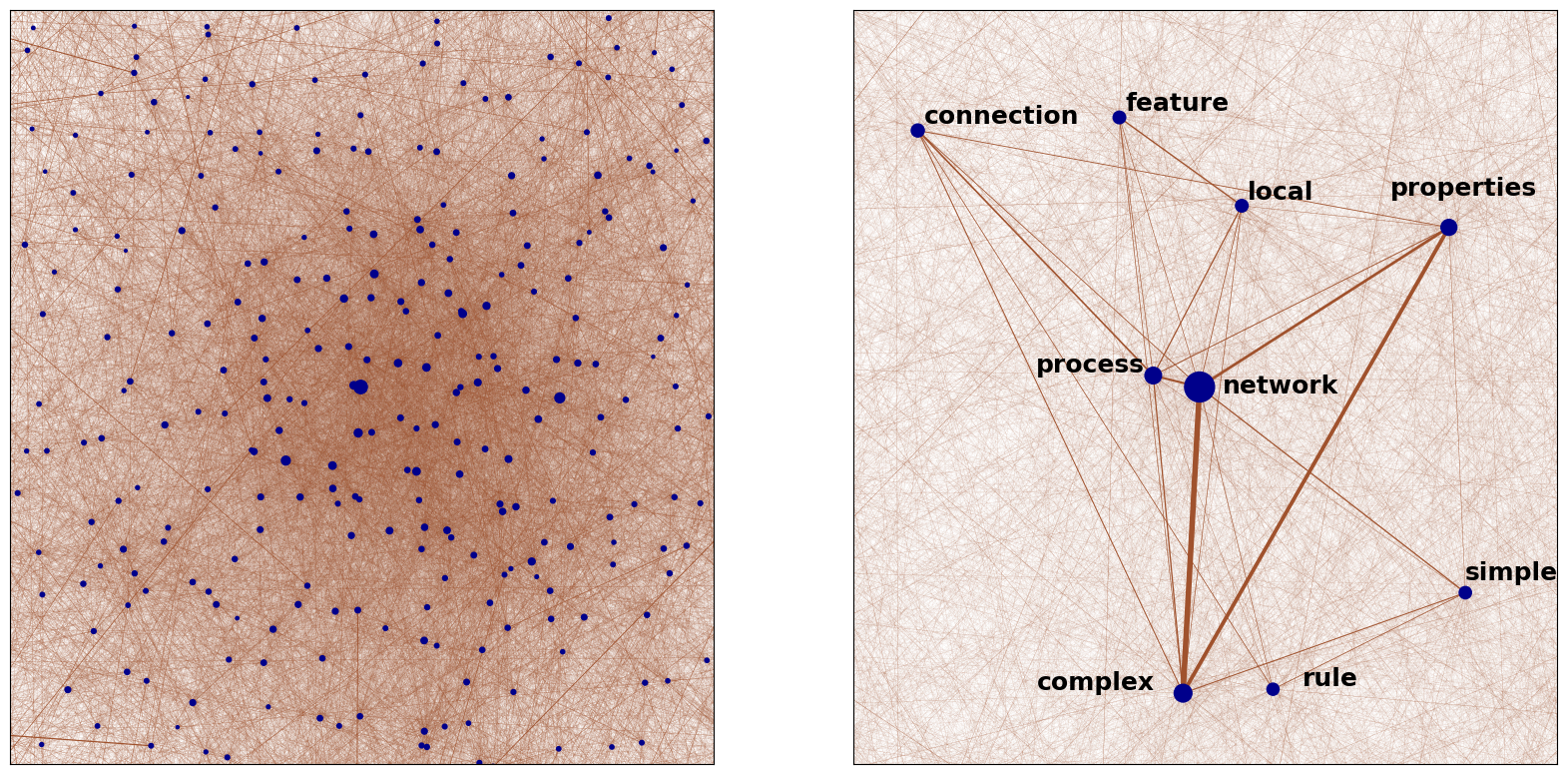}
  \end{center}
 \caption{Homogeneity graph $HG( \mathcal{H})$ of the corpus $\mathcal{H}$ under analysis. On the right side a zoom of the subgraph of neighbors of the word ``network'' is shown.}
\label{Fig6}
 \end{figure}

Figure~\ref{Fig6} shows the homogeneity graph corresponding to the corpus considered, in which the thickness of each edge is proportional to its weight. On the other hand, as it can be seen in the right part of Figure~\ref{Fig6}, there is no link between ``features'' and ``properties'' because $\frac{\partial \mathcal{H}}{\partial \{feature,properties\}}=+\infty $, and the edge joining ``networks'' and ``complex'' is thicker than the rest.

 \begin{figure}[ht]
  \begin{center}
    \includegraphics[width=1\textwidth]{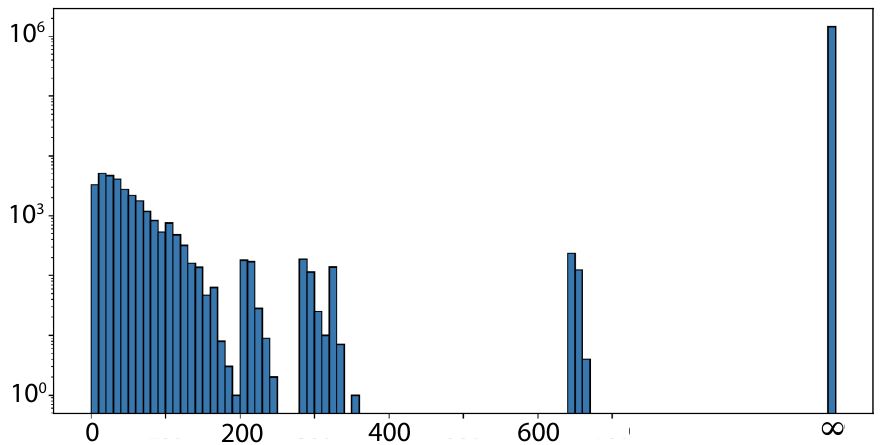}
  \end{center}
 \caption{Histogram clustering the number of pairs of words $\{i,j\}$ by the value of  $\frac{\partial \mathcal{H}}{\partial \{i,j\}}$.}
\label{Fig3}
 \end{figure} 
 
 Also, as it can be seen in the histogram of Figure~\ref{Fig3}, there are more than $10^{3}$ pairs of words $\{i,j\}$ such that $0 \le \frac{\partial \mathcal{H}}{\partial \{i,j\}} \le 10$ and more than $10^{6}$ pairs of words $\{i,j\}$ whose derivative is $+\infty$
 (note that in Figure~\ref{Fig3}, the length of the intervals of the horizontal axis is $10$).

To conclude this section, we would like to point out that the automatic extraction of the linguistic level of a corpus, the search for lexical patterns in sentences of a given author or writer of a particular specialty language, the search for similarities and differences in a set of texts and the automatic classification of texts according to these differences or similarities are some of the possible applications of the methodology underlying this model.




\section{On lexical density and three different rankings of sentences: computational results} \label{sec:experiments}
As far as it is known, the personalized PageRank of a  individual term (node) $i$ is the $i$-component of the stationary state $\pi_0\in\mathbb{R}^n$ ($\|\pi_0\|=1$) of the random walker with transition matrix \cite{Boldi09,brin,page-brin,GaPeRo}
\[
P=\alpha P_B^T+(1-\alpha)\mathbf{v}\mathbf{e}^T,
\]
where $\alpha\in (0,1)$, $B=(b_{ij})$ is the adjacency matrix of the network under consideration, $\mathbf{e}^T=(1,\cdots,1)$, $\mathbf{v}\in\mathbb{R}^n$ ($\|\mathbf{v}\|=1$) is the personalization vector and 
\[
P_B=(p_{ij})=\left(\frac {b_{ij}}{\sum_kb_{ik}}\right).
\]

To carry out our study on the hypergraph $\mathcal{H}$ in which the vertices are the lexical words of the corpus, and the hyperedges are the phrases (sets of lexical words of the corpus located between two periods), we will use the same methodology as in \cite{CABARC20} and \cite{CABARC21} to associate its corresponding PageRank to each node, with the idea of ranking the lexical words according to their importance \cite{Boldi09,brin,page-brin,LangMey06,PRC16}.  For this purpose, taking into account that for the PageRank calculation used throughout this work we have used the algorithm described in {~\cite{ACGPR19}, we will apply this algorithm on three different structures obtained from the application of three different criteria:

\begin{enumerate}
\item $\mathbf{Ranking ~1.}$ To calculate this ranking, we first have built a graph on which to apply the PageRank algorithm. 
In order to do that,  we convert each hyperedge of $\mathcal{H}$ into a clique to obtain  the projection graph $\Pi_2(\mathcal{H})$. After this, taking into account that the average number  of  words of a sentence within the corpus under study  is $5.809$ and that, therefore, the local lexical density is  $5.809$, we can deduce that the damping factor corresponding to this configuration is $0.853$, since 
\[
\begin{split}
5.809= & \mathbb{E}(\ell)=\sum_{k=0}^{\infty}k \cdot \mathbb{P}(\ell=k)
=\sum_{k=1}^{\infty}k \cdot (1-q) \cdot q^{k}\\
=&(1-q)\cdot q\sum_{k=1}^{\infty}k \cdot  q^{k-1}= \frac{q}{1-q}.
\end{split}
\]

\item $ \mathbf{Ranking ~2.}$ To calculate this ranking, we will apply the PageRank algorithm considered on the network  $\Pi_2(\mathcal{H}^*)=L(\mathcal{H})$ so that, once the numerical value attributed to each phrase has been obtained, this value is distributed proportionally among the words that make up that sentence. It is important to note that, in this case, the network considered is a directed network, and that, if $s_1, s_2 \in L(\mathcal{H})$, these sentences will be connected if they have at least one lexical word in common, so that the edge weight $w(s_1\rightarrow s_2)$ is the number of words shared by both sentences multiplied by the number of times that sentence $s_2$ appears repeated in the corpus. Obviously, the edge weight $w(s_1\rightarrow s_2)$  may be different from $w(s_2\rightarrow s_1)$. Now, using the same reasoning as in the previous case, and having in mind that the average number of sentences of a paper included in the corpus under study  is $27.12$, in this context,  the damping factor corresponding to this configuration is $0.96$.

\item $ \mathbf{Ranking ~3.}$ To calculate this ranking, we will apply the PageRank algorithm considered on the weighted graph $HG(\mathcal{H})$. Taking into account that the average number of  words of a sentence is $5.756$ (since, after collapsing words pairs $ \{i,j\}$ such that $\frac{\partial \mathcal{H}}{\partial \{i,j\}}= 0$, the average length of sentences decreases, albeit slightly), the damping factor corresponding to this configuration is $0.852$. Figure~\ref{Fig6} shows the homogeneity graph corresponding to the corpus considered. The size of the nodes is proportional to the component of the PageRank vector corresponding to that node, and the thickness of each edge is proportional to its weight.

\end{enumerate}

 In all the described cases the corresponding value of $q$ is the probability that a random walker will not vary its trajectory by moving to a node directly connected by an edge to the current node instead of jumping to another node in this network not necessarily connected to the previous one. In our situation, this jump can be understood as the end of the current sentence and the starting point of a new sentence for Ranking 1 and Ranking 3, and  as the end of the current paper and the starting point of a new paper for Ranking 2. To complete the necessary elements to apply the algorithm, we will point out that for Ranking 1 and Ranking 3 the personalization vector considered is the (relative) frequency of lexical words, and for Ranking 2 the personalization vector considered is the (relative) frequency of each sentence  included in the corpus under study.
 

\begin{table}[htbp]
\begin{center}
\resizebox{6cm}{!} {
\begin{tabular}{|c|cccc|}
\hline 
&  Ranking 1 & Ranking 2&   Ranking 3&\\ \hline 
1st & network & network &network&\\ \hline 
2nd & system  & system& system& \\ \hline
3rd &  model & model & model &\\ \hline
4th & complex & complex & complex&\\ \hline
5th & process& number & graph&\\ \hline
6th &  number & process & process&\\ \hline
7th &  information  & structure & structure&\\ \hline
8th & graph & new&information& \\ \hline
9th &  new& information &number&\\ \hline
10th & structure & distribution &new& \\ \hline
11th &  properties & properties &properties&\\ \hline
12th &  distribution  & graph &distribution &\\ \hline
13th & study & study&dynamics &\\ \hline
14th &  dynamics& dynamics &study &\\ \hline
15th & case & interaction & analysis& \\ \hline
\end{tabular}
}
\caption{Rankings of lexical words}
\label{Table1}
\end{center}
\end{table}

 \begin{figure}[h!]
  \begin{center}
    \includegraphics[width=0.8\textwidth]{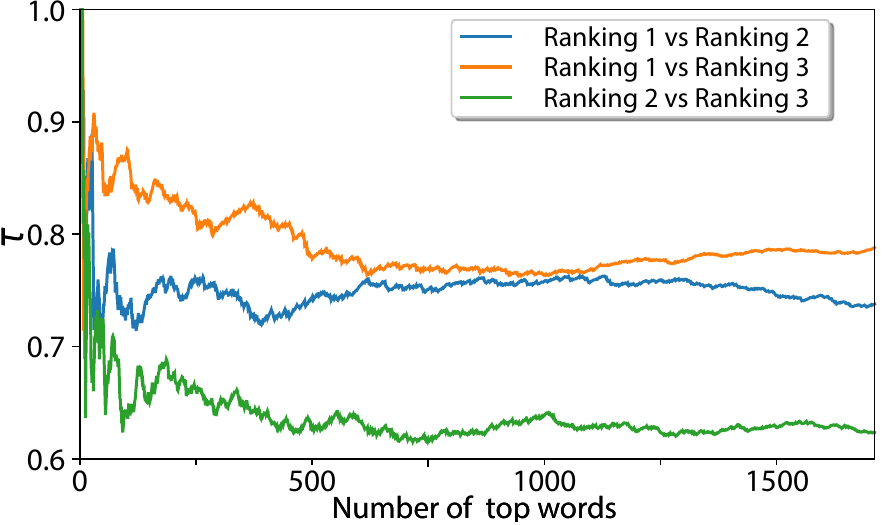}
  \end{center}
 \caption{Kendall's tau coefficient variation by comparing the rankings in pairs among them, depending on the number of top lexical words considered in each ranking.}
\label{Fig9}
 \end{figure}

As it can be seen in Table~\ref{Table1}, there is hardly any difference at the top of the three rankings. As expected, Ranking 3 gives us the most representative words of the corpus in the sense that they are the words at the heart of the corpus linking the largest number of sentences together. In any case, the three rankings should not be very different from each other, as it is actually the case (since the first four positions are occupied by the same words in all three cases) and, as it happens in the case under study, Ranking ~1 and Ranking ~3 are more similar to each other than to Ranking ~2. However, as the number of words considered at the top of each ranking increases, the differences between the three rankings become much more evident, as it can be seen in Figure~\ref{Fig9}, where we plot the differences between these rankings by visualizing the variation of the Kendall's tau coefficient ($\tau$) \cite{ken} regarding the number of lexical words considered in the three rankings.

\section{Seeking for distinctive characteristics that allow distinguishing the styles of different authors and language levels}\label{sec:seeking}
By considering several types and models of hypergraphs and hyperstructures for a given text or corpus, we can associate to that written text or corpus various features that allow us to identify it as if it were some sort of mathematical signature associated with them. For example, for a given text it is possible to consider a hypergraph in which the nodes are the words and the hyperedges are the sentences, another in which the nodes are the words and the hyperedges are the paragraphs, another in which the nodes are the sentences and the hyperedges are the paragraphs, just to mention some of the possibilities.
This succession of mathematical structures and the different parameters (such as diameter, degree distribution, centrality, efficiency, among others, that characterize them) are, without a doubt, elements that configure and allow us to characterize and compare different texts, making it clear the characteristics that constitute their seal of identity in terms of style.

\section{Conclusions}\label{sec:conclu}

We introduce and study the derivative of a hypergraph and the homogeneity graph of a hypergraph as new and useful structures that can be used to study the degree of independence of the nodes of a hypergraph as well as to obtain a ranking of the most representative nodes of the hypergraph in the sense that the lexical words represented by these nodes link the most significant ideas and concepts of the text without necessarily being those terms usually considered as keywords.

These concepts allow us to associate not only a numerical index that allows us to quantify the heterogeneity and lack of similarity between the nodes of the hypergraph, but also to associate a graph to characterize the heterogeneity and dissimilarity of the different elements of the considered hypergraph. 

Moreover, these concepts also allow us to obtain technical characteristics related to the styles of the different authors and the language competence level of any text written in English, as well as their possible application to text classification, text summarization, automated translation, stylometry and authorship detection.  

Undoubtedly, the tools derived from the linguistic analysis obtained by using this new tool will provide with new models and better instruments to typify and locate the characteristics of the style of different authors together with the style and intrinsic linguistic characteristics found in specialized texts in terms of collocations, word sense desambiguation and syntagmatic structures. 

Finally, it is important to mention that the construction of tools to find lexical patterns of the style of an author or a text belonging to a specialty language, the automatic classification of texts according to their style and the automatic labeling and identification/verification of lexical patterns are some possible additional applications of these new tools. 
\section*{Acknowledgements}

Authors would like to thanks Karin Alfaro-Bittner for some inspiring discussions. This work has been partially supported by projects PGC2018-101625-B-I00 (Spanish Ministry, AEI/FEDER, UE) and M1993 Grant (Rey Juan Carlos University, Spain). Authors acknowledge the usage of the resources, technical expertise and assistance provided by the supercomputer facility CRESCO of ENEA in Portici (Italy).


\end{document}